\documentclass[11pt]{article}
\usepackage[margin=1in]{geometry}
\usepackage{amsfonts,amsmath,amssymb}
\usepackage{graphbox}

\PassOptionsToPackage{obeyspaces}{url}
\usepackage[colorlinks=true,citecolor=blue,urlcolor=blue,linkcolor=blue,bookmarksopen=true]{hyperref}
\usepackage{amsthm}
\usepackage{tikz}
\usepackage{tikz-layers}
\usetikzlibrary{calc,math,positioning,shapes,decorations.pathmorphing}

\usepackage{etoolbox}
\patchcmd{\thebibliography}{\leftmargin\labelwidth}{\leftmargin\labelwidth\addtolength\itemsep{-0.1\baselineskip}}{}{}

\author{Joseph Briggs\thanks{Dept.\ of Mathematics \& Statistics, Auburn University, Auburn, AL, USA. \texttt{\{jgb0059,coc0014\}@auburn.edu}.} \and Chris Wells\footnotemark[1]}

\title{Phase transitions in isoperimetric problems on the integers}
\date{}

\usepackage[nameinlink,sort]{cleveref}

\newtheorem{theorem}{Theorem}
\newtheorem{lemma}[theorem]{Lemma}

\crefname{conj}{conjecture}{conjectures}

\crefname{claim}{claim}{claims}

\newtheorem{prop}[theorem]{Proposition}
\crefname{prop}{proposition}{propositions}

\theoremstyle{definition}

\crefname{defn}{definition}{definitions}

\crefname{remark}{remark}{remarks}

\newtheorem{question}[theorem]{Question}
\crefname{question}{question}{questions}

\crefname{enumi}{part}{parts}

\numberwithin{theorem}{section}

\makeatletter
\DeclareRobustCommand{\crefnosort}[1]{%
    \begingroup\@cref@sortfalse\cref{#1}\endgroup
}
\DeclareRobustCommand{\Crefnosort}[1]{%
    \begingroup\@cref@sortfalse\Cref{#1}\endgroup
}
\makeatother

\newcommand*{\eqdef}{\stackrel{\mbox{\normalfont\tiny{def}}}{=}}        
\newcommand*{\abs}[1]{\lvert #1\rvert}                
\newcommand*{\abss}[1]{\bigl\lvert #1\bigr\rvert}     

\renewcommand*{\epsilon}{\varepsilon}       
\newcommand*{\Z}{\mathbb{Z}}                
\newcommand*{\N}{\mathbb{N}}                
\newcommand*{\e}{\mathbf{e}}
\DeclareMathOperator{\cay}{Cay}

\newcommand*{\boundary}[2]{\operatorname{\partial_{#1,{\mathnormal #2}}}}
\newcommand*{\opt}[2]{\operatorname{\partial^*_{#1,{\mathnormal #2}}}}
\newcommand*{\Opt}[2]{\operatorname{Opt_{#1,{\mathnormal #2}}}}

\begin{document}
\maketitle

\begin{abstract}
    Barber and Erde asked the following question: if $B$ generates $\Z^d$ as an additive group, then must the extremal sets for the vertex/edge-isoperimetric inequality on the Cayley graph $\cay(\Z^d,B)$ form a nested family?
    We answer this question negatively for both the vertex- and edge-isoperimetric inequalities, specifically in the case of $d=1$.
    The key is to show that the structure of the cylinder $\Z\times(\Z/k\Z)$ can be mimicked in certain Cayley graphs on $\Z$, leading to a phase transition.
    We do, however, show that Barber--Erde's question for Cayley graphs on $\Z$ has a positive answer if one is allowed to ignore finitely many sets.
\end{abstract}

\section{Introduction}

Let $G$ be a directed graph and $A\subseteq V(G)$.
The \emph{edge-boundary} of $A$ is the set
\[
    \boundary e G A\eqdef\{(u,v)\in E(G): u\in A, v\notin A\},
\]
and the \emph{vertex-boundary} of $A$ is the set
\[
    \boundary v G A\eqdef\{v\in V(G)\setminus A: \text{there is }u\in A\text{ with }(u,v)\in E(G)\}.
\]
These definitions naturally specialize to the case when $G$ is an undirected graph by imagining that each edge $\{u,v\}$ is formed by the two directed edges $(u,v)$ and $(v,u)$.
Throughout this paper, every graph is assumed to be directed.
\medskip

An \emph{isoperimetric problem} on the graph $G$ is concerned with determining the minimum value of $\abs{\boundary eGA}$ or $\abs{\boundary vGA}$ over all subsets $A\subseteq V(G)$ of a fixed size.
To this end, define
\begin{align*}
    \opt e G (n) &\eqdef \min\bigl\{\abs{\boundary eGA}:A\subseteq V(G),\ \abs{A}=n\bigr\},\qquad\text{and}\\
    \opt v G(n) &\eqdef \min\bigl\{\abs{\boundary vGA}:A\subseteq V(G),\ \abs{A}=n\bigr\}.
\end{align*}
In particular, for every finite $A\subseteq V(G)$, we have $\abs{\boundary e G A}\geq\opt{e}{G}(\abs A)$ and $\abs{\boundary{v}{G}A}\geq\opt{v}{G}(\abs A)$.
\medskip

This paper is concerned with isoperimetric problems in Cayley graphs on integer lattices.
For an additive group $\Gamma$ and a subset $B\subseteq\Gamma$, the Cayley graph $\cay(\Gamma,B)$ is the graph with vertex-set $\Gamma$ where $(u,v)$ is an edge if and only if $v-u\in B$.
Note that we may consider $\cay(\Gamma,B)$ to be an undirected graph if and only if $B=-B$; such a set $B$ is said to be \emph{symmetric}.

A subset $B\subseteq\Gamma$ is said to \emph{generate} $\Gamma$ if every element of $\Gamma$ can be written as an integer combination of elements of $B$.
Equivalently, $B$ generates $\Gamma$ if the Cayley graph $\cay(\Gamma,B)$ is weakly connected.
\medskip

Motivated by known isoperimetric inequalities on the integer lattice, Barber and Erde asked the following question:
\begin{question}[Question 1 of Barber--Erde~\cite{barber_iso}]\label[question]{q.main}
    Let $B\subseteq\Z^d$ be a finite set which generates $\Z^d$ and set $G=\cay(\Z^d,B)$.
    \begin{itemize}
        \item Is there always an ordering $v_1,v_2,\ldots$ of $\Z^d$ for which $\abs{\boundary{e}{G}\{v_1,\dots,v_n\}}=\opt{e}{G}(n)$ for every $n\in\N$?
        \item Is there always an ordering $u_1,u_2,\ldots$ of $\Z^d$ for which $\abs{\boundary{v}{G}\{u_1,\dots,u_n\}}=\opt{v}{G}(n)$ for every $n\in\N$?
    \end{itemize}
\end{question}

The classical result of Bollob\'as--Leader~\cite{bollobas_edge} answers this affirmatively for $\opt{e}{G}(n)$ when $G$ is the $\ell_1$-lattice on $\Z^d$.
The two celebrated orderings of Wang--Wang~\cite{wang_l1} and Radcliffe--Veomett~\cite{radcliffe_vertex} likewise answer this affirmatively for $\opt{v}{G}(n)$ when $G$ is the $\ell_1$- and the $\ell_\infty$-lattice on $\Z^d$, respectively.
Beyond this, Gupta--Levcovitz--Margolis--Stark~\cite{gupta_planar} completely classify the optimizers for the vertex-isoperimetric inequality on the $\ell_1$-lattice on $\Z^2$.

In addition to these known orderings, sets of lattice points witnessing $\opt{v}{G}(n)$ and $\opt{e}{G}(n)$ \emph{asymptotically} for any Cayley graph $G$ on $\Z^d$ were established by Ruzsa~\cite{ruzsa_asymptotics} and Barber--Erde~\cite{barber_iso}, respectively; the latter confirmed a conjecture of Tsukerman--Veomett~\cite{tsukerman_conjecture}.
A related general lower bound for Cayley graphs on \emph{finite} abelian groups has been established by Lev~\cite{lev_edge}.
It is this evidence that led Barber and Erde to ask \Cref{q.main}.

Since Barber and Erde asked \Cref{q.main}, Barber--Erde--Keevash--Roberts~\cite{barber_stability} have shown that these asymptotic optimizers are \emph{stable}; that is, for each $i \in \{\operatorname{e}, \operatorname{v}\}$, any two sets $A,A'$ of size $n$ with $\boundary{\mathnormal i}{G} A,\boundary{\mathnormal i}{G} A'\leq\bigl(1+o(1)\bigr)\opt{\mathnormal i}{G}(n)$ are at an edit distance of $o(n)$ from one another.
This astonishing result yields extra evidence to support \Cref{q.main}.

However, surprisingly, we show that the answer to \Cref{q.main} is ``no'' already in the simplest case of $d=1$.

For $i\in\{\operatorname e,\operatorname v\}$, let $\Opt{{\mathnormal i}}{G}(n)$ denote the set of all $A\subseteq V(G)$ with $\abs A=n$ and $\abs{\boundary{\mathnormal i}{G}A}=\opt{\mathnormal i}{G}(n)$.
A restatement of \Cref{q.main} is: For $i\in\{\operatorname e,\operatorname v\}$, can we find $A_n\in\Opt{\mathnormal i}{G}(n)$ such that $A_1\subseteq A_2\subseteq A_3\subseteq\cdots$?
\begin{theorem}\label[theorem]{not_nested}
    For every integer $N$ and for each $i\in\{\operatorname e,\operatorname v\}$, there is a subset $B\subseteq\Z$ and integers $n_1,n_2$ satisfying
    \begin{itemize}
        \item $N\leq n_1\leq n_2-N$, and
        \item $B$ is finite, symmetric and generates $\Z$, and
        \item With $G=\cay(\Z,B)$, no element of $\Opt{\mathnormal i}{G}(n_1)$ is a subset of any element of $\Opt{\mathnormal i}{G}(n_2)$.
    \end{itemize}
\end{theorem}

Despite this, and less surprisingly, we show that \Cref{q.main} has a positive answer \emph{eventually} for Cayley graphs on $\Z$.
This answers Question 2 of Barber--Erde~\cite{barber_iso} positively in the case of $d=1$.

\begin{theorem}\label[theorem]{interval_opt}
    Suppose that $B\subseteq\Z$ is any finite set which generates $\Z$ and set $G=\cay(\Z,B)$.
    There is an integer $N=N(B)$ such that $[n]\in\Opt{e}{G}(n)$ and $[n]\in\Opt{v}{G}(n)$ whenever $n\geq N$.
\end{theorem}

We note that, since an earlier version of this paper was written, a family of examples in dimensions $d \geq 2$ having non-nestability of edge-isoperimetric optimizers has been constructed by Strachan and Swanepoel in \cite{strachan2025edge}. Taking $B=\{\pm e_1, \pm e_2, \dots, \pm e_d, \pm 2e_1\} \subset \Z^d$, they showed that an analogue of \Cref{interval_opt} does not hold for $d \geq 2$ for $i=e$ (but the $i=v$ case is still open, see \Cref{sec:remarks}). Although not explicitly written there, their construction also gives rise to arbitrarily wide gaps between non-nested optimizers, in the sense of \Cref{not_nested}.

This paper is organized as follows.
The entirety of \Cref{sec:construction} is dedicated to proving \Cref{not_nested}.
\Cref{not_nested} is explicitly proved in \Cref{sec:not_nested}, but the majority of the details involved are established in \Cref{sec:edge,sec:vertex}.
After detailing the constructions, we prove \Cref{interval_opt} in \Cref{sec:interval_opt}.
We conclude by discussing further questions in \Cref{sec:remarks}.

\subsection{Notation and preliminaries}
Let $\Gamma$ be an additive group.
For sets $A,B\subseteq\Gamma$, the Minkowski sum and difference are $A+B\eqdef\{a+b:a\in A,\ b\in B\}$ and $A-B\eqdef\{a-b:a\in A,\ b\in B\}$, respectively.
If $B=\{b\}$, then we instead write $A+b=A+\{b\}$ for simplicity.
We additionally abbreviate $A\pm B=A+(B\cup(-B))$.
Combining these two abbreviations, $A\pm b=A+\{b,-b\}$.
Finally, for $k\in\Z$, the $k$-dilate of $A$ is $kA\eqdef\{ka:a\in A\}$.\footnote{In some additive combinatorics texts, $kA$ is defined as the $k$-fold sum of $A$ with itself. Since we have no need of this operation, we instead use $kA$ as the $k$-dilate to imitate the natural notation of $k\Z$ to denote the integers divisible by $k$.}

In the Cayley graph $G=\cay(\Gamma,B)$, observe that $\boundary{v}{G}A=(A+B)\setminus A$ and that $\boundary{e}{G}A$ can be naturally identified with the set $\bigcup_{b\in B}\bigl(\{b\}\times\bigl((A+b)\setminus A\bigr)\bigr)$.
\medskip

For integers $m,n$, we write $[m,n]\eqdef\{m,m+1,\dots,n-1,n\}$ to be the integers between $m$ and $n$, inclusive.
If $m>n$, then $[m,n]=\varnothing$.
As is common, we abbreviate $[n]\eqdef[1,n]$.
\medskip

An \emph{arithmetic progression} with step size $k$ is any set of the form $\{x,x+k,\dots,x+nk\}$ for some $x\in\Z$ and $n\in\Z_{\geq 0}$.

\section{Construction of non-nested extremal sets}\label{sec:construction}
In order to build counterexamples to Barber--Erde's question, we rely on the known isoperimetric inequalities for the $\ell_1$- and $\ell_\infty$-grids on $\Z^2$.
Set
\[
    \ell_1^2 \eqdef\cay(\Z^2,\pm\{\e_1,\e_2\}),\qquad\text{and}\qquad
    \ell_\infty^2 \eqdef\cay(\Z^2,\{0,\pm 1\}^2 \setminus \{\mathbf 0\}).
\]
Here $\e_1=(1,0)$ and $\e_2=(0,1)$.
In particular, for $p\in\{1,\infty\}$, the graph $\ell_p^2$ has vertex set $\Z^2$ where $(\mathbf x,\mathbf y)$ is an edge if and only if $\lVert \mathbf x-\mathbf y\rVert_p=1$.

\begin{theorem}[Bollob\'as--Leader~\cite{bollobas_edge}]\label[theorem]{1_iso}
    If $A\subseteq\Z^2$ is any set of size $n$, then $\abs{\boundary{e}{\ell_1^2} A}\geq 4\sqrt n$ with equality if and only if $n=k^2$ and $A$ is a translate of $[k]^2$ for some integer $k$.
\end{theorem}
Bollob\'as and Leader do not explicitly state this fact in their paper, but it can be quickly deduced from the second paragraph in the proof of Theorem 15 in \cite{bollobas_edge}.

\begin{theorem}\label[theorem]{infty_iso}
    If $A\subseteq\Z^2$ is any set of size $n$, then $\abs{\boundary{v}{\ell_\infty^2} A}\geq 4(\sqrt n+1)$ with equality if and only if $n=k^2$ and $A$ is a translate of $[k]^2$ for some integer $k$.
\end{theorem}
Radcliffe and Veomett~\cite{radcliffe_vertex} proved a generalization of the inequality in the above theorem to arbitrary dimensions, but they did not establish the uniqueness of the extremal sets when $n=k^2$.
Since we will need to know that translates of $[k]^2$ are the \emph{unique} extremal sets, we give an explicit proof of \Cref{infty_iso} in \Cref{sec:isoproofs}.

\subsection{Mimicking the cylinder}\label{sec:motivation}
Consider the $\ell_1$- and $\ell_\infty$-cylindrical grids:
\[
    C_1(b) \eqdef\cay(\Z\times(\Z/b\Z),\pm\{\e_1,\e_2\}),\qquad\text{and}\qquad
    C_\infty(b) \eqdef\cay(\Z\times(\Z/b\Z),\{0,\pm 1\}^2\setminus\{\mathbf 0\}).
\]
While we are unaware of an explicit result in this direction, the edge-isoperimetric problem on $C_1(b)$ and the vertex-isoperimetric problem on $C_\infty(b)$ each exhibit a ``phase transition'', meaning that the extremal sets are not nested.
Indeed, if $\sqrt n\ll b$, then the edge-isoperimetric problem on $C_1(b)$ (resp.\ the vertex-isoperimetric problem on $C_\infty(b)$) should be very similar to that for $\ell_1^2$ (resp.\ $\ell_\infty^2$) in that the extremal examples should be translates of the square $[\sqrt n]^2$.
On the other hand, if $\sqrt n\gg b$, then the extremal examples for the edge-isoperimetric inequality on $C_1(b)$ (resp.\ the vertex-isoperimetric inequality on $C_\infty(b)$) should be translates of the ``interval'' $[n/b]\times(\Z/b\Z)$.

\begin{figure}
    \includegraphics[align=c,width=0.47\textwidth]{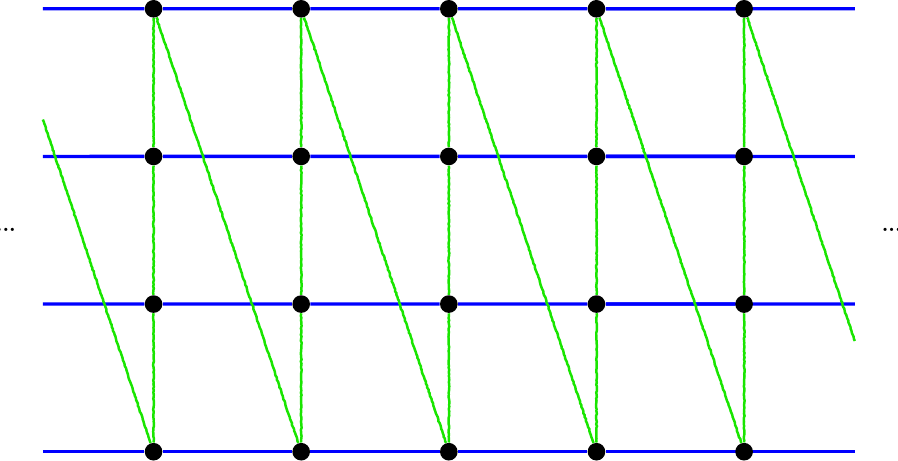}
    \hfill
    \includegraphics[align=c,width=0.47\textwidth]{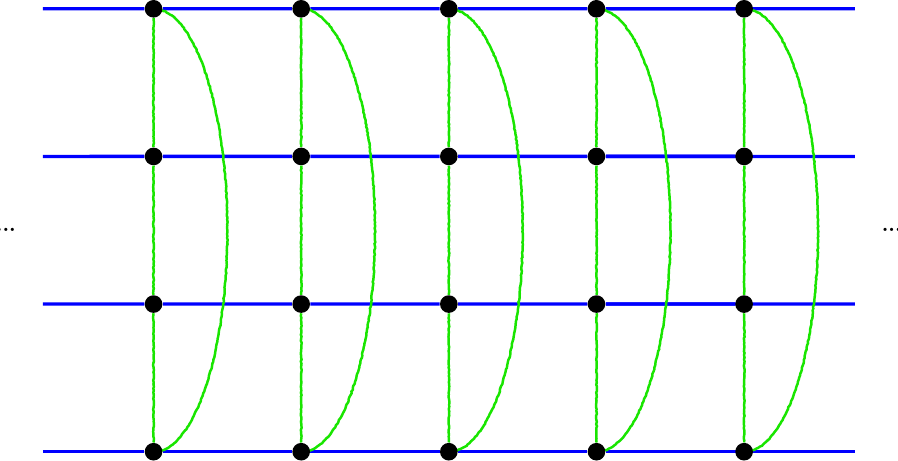}
    \caption{\label{fig:l1cylinder}The graphs $\cay(\Z,\pm\{1,4\})$ (left) and $C_1(4)$ (right).}
\end{figure}

\begin{figure}
    \includegraphics[align=c,width=0.47\textwidth]{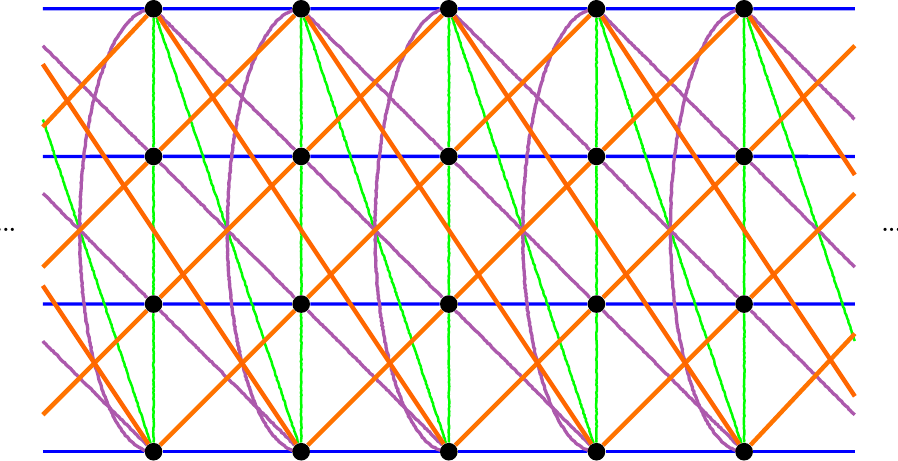}
    \hfill
    \includegraphics[align=c,width=0.47\textwidth]{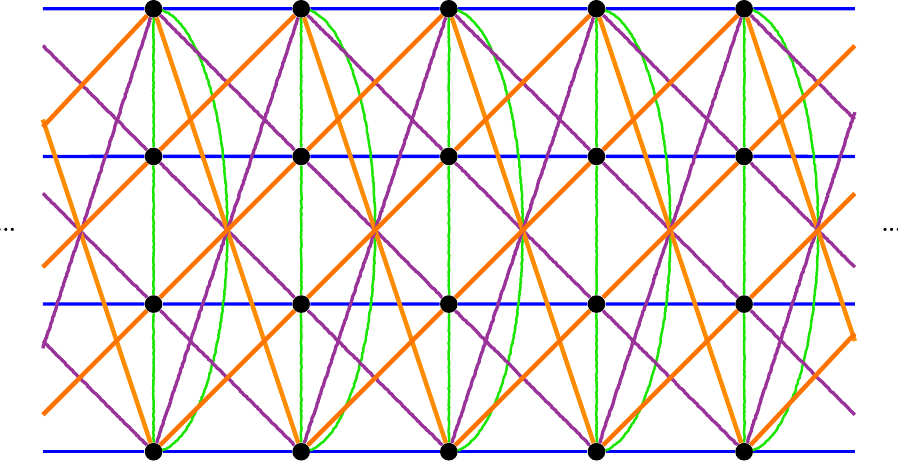}
    \caption{\label{fig:linftycylinder}The graphs $\cay(\Z,\pm\{1,3,4,5\})$ (left) and $C_\infty(4)$ (right).}
\end{figure}

Now, consider the map $\phi\colon\Z\to\Z\times(\Z/b\Z)$ defined by $\phi(x)=(q,r)$ where $x=qb+r$ is the Euclidean decomposition of $x$.
Notice that, $\phi(b)=\e_1$ and $\phi(1)=\e_2$.
Due to this, $\phi$ is a ``not-quite-isomorphism'' between $\cay(\Z,\pm\{1,b\})$ and $C_1(b)$ (see \Cref{fig:l1cylinder}).
Therefore, we might expect that the edge-isoperimetric problem on $\cay(\Z,\pm\{1,b\})$ behaves similarly to that on $C_1(b)$.
Similarly, but perhaps more subtly, $\phi$ is also a ``not-quite-isomorphism'' between $\cay(\Z,\pm\{1,b-1,b,b+1\})$ and $C_\infty(b)$ (see \Cref{fig:linftycylinder}) and so we might expect that the vertex-isoperimetric inequalities on these two graphs behave similarly.

We do not explicitly make use of these cylindrical grids in our arguments, but it will be helpful to keep this motivation in mind.
While there is a reasonable way to directly link the edge-isoperimetric problems on $\cay(\Z,\pm\{1,b\})$ and $C_1(b)$, doing the same for the vertex-isoperimetric problems on $\cay(\Z,\pm\{1,b-1,b,b+1\})$ and $C_\infty(b)$ is more complicated.
Instead, the proofs that follow imitate the methods that one would likely use in order to show that, e.g., the vertex-isoperimetric problem on $C_\infty(b)$ and $\ell_\infty^2$ align when $\sqrt n\ll b$.

We remark that that our arguments can be phrased in terms of (a variation on) Freiman homomorphisms from additive combinatorics (see, e.g., {\cite[Chapter 7]{zhao_book}}).
However, since our arguments are so rudimentary, we feel that phrasing them in this language would only introduce unnecessary complications.

\subsection{Construction for the edge-isoperimetric problem}\label{sec:edge}
Following the motivation set out in \Cref{sec:motivation}, our goal is to relate the edge-isoperimetric problems on $\cay(\Z,\pm\{1,b\})$ and $\ell_1^2$.
The key idea is to show that if $A\subseteq\Z$ is optimal, then either $A$ is an interval or we can use the Euclidean map $x\mapsto(q,r)$ where $x=qb+r$ to embed $A$ into $\ell_1^2$ without affecting its boundary.
To accomplish the latter, we will need to ``cut'' along a residue class modulo $b$.
The following proposition contains the main observation which enables this feat: if an optimal $A$ is not an interval, then $A$ must be missing some residue modulo $b$.

\begin{prop}\label[prop]{every_res_then_interval}
    Fix finite $A,B\subseteq\Z$ and set $G=\cay(\Z,B)$.
    If $A$ contains every residue modulo $\abs b$ for each $b\in B$, then $\abs{\boundary{e}{G} A}\geq \sum_{b\in B} \abs b$.
    Furthermore:
    \begin{enumerate}
        \item Equality holds if $A$ is an interval.
        \item If $B$ contains either $1$ or $-1$, then equality holds if and only if $A$ is an interval.
    \end{enumerate}
\end{prop}
\begin{proof}
    For $b\in B$ and $r\in[0,\abs b-1]$, set $A_{b,r}=A\cap(r+\abs b\Z)$.
    Of course, if $a\in A_{b,r}$ and $a+b\in A$, then $a+b\in A_{b,r}$ as well.
    Therefore,
    \[
        \abs{\boundary{e}{G}A}=\sum_{b\in B}\abs{(A+b)\setminus A}=\sum_{b\in B}\sum_{r=0}^{\abs b-1}\abs{(A_{b,r}+b)\setminus A_{b,r}}.
    \]
    Now, since $A_{b,r}$ is non-empty and finite for each $b,r$ by assumption, it is quick to observe that $\abs{(A_{b,r}+b)\setminus A_{b,r}}\geq 1$ with equality if and only if $A_{b,r}$ is an arithmetic progression with step size $\abs b$.
    This concludes the proof.
\end{proof}

Already this is enough to demonstrate the phase transition in the edge-isoperimetric problem on $\cay(\Z,\pm\{1,b\})$.
As was mentioned above, the key idea is to take an optimal $A$, which is not an interval, and ``cut'' along some missing residue class modulo $b$ in order to embed $A$ into $\ell_1^2$ without affecting its boundary.

\begin{figure}
    \includegraphics[align=c,width=0.45\textwidth]{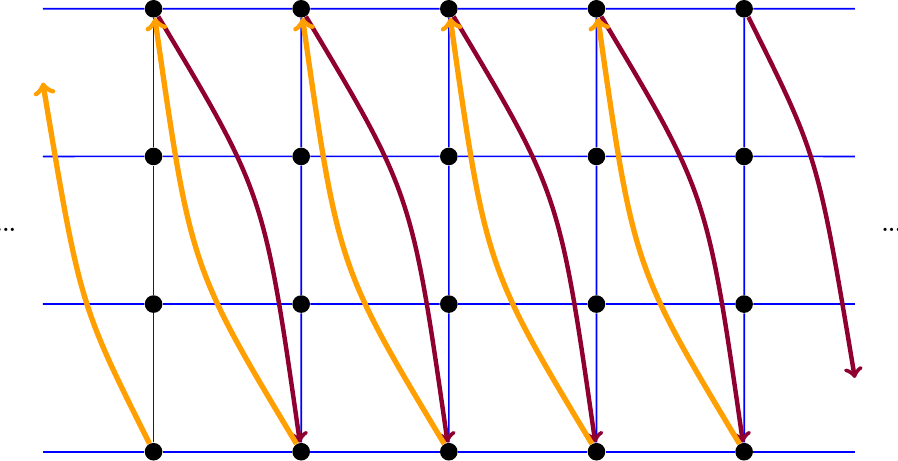}
    \hfill$\stackrel{\mbox{\normalfont\small{$\phi'$}}}{\longrightarrow}$\hfill
    \includegraphics[align=c,width=0.45\textwidth]{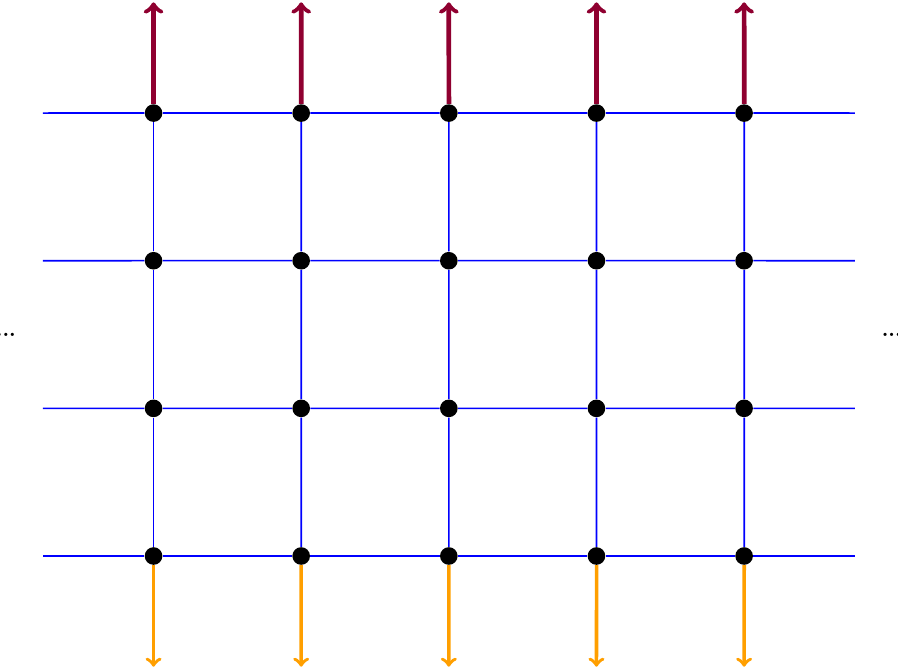}
    \caption{\label{fig:edgemorph}The map $\phi'$ from the proof of \Cref{edge_not_nested} with $b=4$.}
\end{figure}

\begin{theorem}\label[theorem]{edge_not_nested}
    Fix an integer $b\geq 2$ and set $G=\cay(\Z,\pm\{1,b\})$.
    For every positive integer $n$:
    \begin{enumerate}
        \item If $\sqrt n>{b+1\over 2}$, then every element of $\Opt{e}{G}(n)$ is an interval of length $n$.
        \item If $n=k^2$ and $k<{b+1\over 2}$, then every element of $\Opt{e}{G}(n)$ is a translate of $[k]+b[k]$.
    \end{enumerate}
\end{theorem}
\begin{proof}
    Observe that
    \begin{equation}\label{eqn:edge_interval}
        \abs{\boundary{e}{G}[n]}=2\bigl(\min\{b,n\}+1\bigr).
    \end{equation}

    Suppose that $A\in\Opt{e}{G}(n)$ is \emph{not} an interval.
    We will show the following:
    \begin{enumerate}
        \item $\sqrt n\leq{b+1\over 2}$, and
        \item If $n=k^2$ and $2\leq k<{b+1\over 2}$, then $\abs{\boundary{e}{G}A}<\abs{\boundary{e}{G}[n]}$ and $A$ is a translate of $[k]+b[k]$.
    \end{enumerate}
    This will establish the claim.
    \medskip

    Thanks to \Cref{every_res_then_interval}, since $A$ is not an interval, we know that $A$ must be void of some residue modulo $b$.
    Since this problem is translation invariant, we may suppose, without loss of generality, that no element of $A$ is divisible by $b$.
    \medskip

    Consider the map $\phi\colon\Z\to\Z^2$ defined by $\phi(x)=(q,r)$ where $x=qb+r$ is the Euclidean decomposition of $x$; note that the image of $\phi$ is $\Z\times[0,b-1]$.
    Since $A$ has no elements which are divisible by $b$, we see that $\phi(A)\subseteq\Z\times[1,b-1]$.
    Now, consider a related map $\phi'\colon E(G)\to E(\ell_1^2)$ defined by
    \begin{align*}
        \phi'(x,x\pm b) &=(\phi(x),\phi(x)\pm\e_1),\quad\text{and}\\
        \phi'(x,x\pm 1) &=(\phi(x),\phi(x)\pm\e_2).
    \end{align*}
    A pictorial representation of $\phi'$ can be seen in \Cref{fig:edgemorph}.
    Note that $\phi'$ does not respect the symmetric nature of these two Cayley graphs since, e.g., $\phi'(0,-1)=\bigl((0,0),(0,-1)\bigr)$ whereas $\phi'(-1,0)=\bigl((-1,b-1),(-1,b)\bigr)$.
    However, since $A$ has no element which is divisible by $b$, we find that $e\in E(G[A])$ if and only if $\phi'(e)\in E(\ell_1^2[\phi(A)])$.
    Therefore, $\phi'\bigl(\boundary e G A\bigr)=\boundary e{\ell_1^2}\phi(A)$.
    In particular, due to the isoperimetric inequality for $\ell_1^2$ (\Cref{1_iso}), we have $\abs{\boundary{e}{G}A}\geq 4\sqrt n$ with equality if and only if $n=k^2$ and $\phi(A)$ is a translate of $[k]^2$.

    Due to \cref{eqn:edge_interval}, this then implies that $4\sqrt n\leq 2(b+1)\implies \sqrt n\leq {b+1\over 2}$, which establishes the first part of the claim.

    Now, suppose that $n=k^2$ and that $k=\sqrt n<{b+1\over 2}$.
    In particular, $k<b$ and so $[k]^2\subseteq\Z\times[1,b-1]$.
    Thus, from the observation above, we must have $\abs{\boundary{e}{G}A}=4k$ with equality if and only if $\phi(A)$ is a translate of $[k]^2$.
    Of course, $\phi(A)$ is a translate of $[k]^2$ if and only if $A$ is a translate of $[k]+b[k]$.
    To conclude, we check that $[n]$ cannot be optimal in this situation, unless $k=1$ in which case $[k]+b[k]=\{1+b\}$.
    To see this, note that $2k<k^2+1$ whenever $k\geq 2$ and so $\abs{\boundary{e}{G}[n]}=2\bigl(\min\{n,b\}+1\bigr)>4k=\abs{\boundary{e}{G}A}$.
\end{proof}

\subsection{Construction for the vertex-isoperimetric problem}\label{sec:vertex}
We now turn our attention to the vertex-isoperimetric problem.

Recalling the motivation from \Cref{sec:motivation}, our goal is to relate the vertex-isoperimtric problems on $\cay(\Z,\pm\{1,b-1,b,b+1\})$ and $\ell_\infty^2$.
Just like in the previous section, the key is to show that if $A\subseteq\Z$ is optimal, then either $A$ is an interval or we can use the Euclidean map $x\mapsto(q,r)$ where $x=qb+r$ to embed $A$ into $\ell_\infty^2$ without affecting its boundary.
We will again look to ``cut'' along a residue class modulo $b$ in order to accomplish this feat.

In the previous section, it was enough to avoid a single residue class modulo $b$, but here we will need more.
When it came to edge-boundaries, we could get away with ``cutting'' along a single residue class since, e.g., the edges $(-1,0)$ and $(1,0)$ are distinct even though their end-points are not; however, when considering only vertex-boundaries, the fact that the end-points are not distinct causes issues.
The following two lemmas, \Cref{missing_none,missing_one}, slowly build up to \Cref{two_consecutive}, which states that if an optimal $A$ is not an interval, then $A$ must be missing \emph{two consecutive} residue classes modulo $b$.
This will then allow us to perform the desired cut.

\begin{lemma}\label[lemma]{missing_none}
    Fix a finite set $B\subseteq\Z_{\geq 1}$, set $b=\max B$ and suppose that $1,b-1\in B$ as well.
    If $A\subseteq\Z$ is a finite set which contains every residue modulo $b$, then $\abs{(A\pm B)\setminus A}\geq 2b$ with equality if and only if $A$ is an interval.
\end{lemma}
\begin{proof}
    For each $r\in[0,b-1]$, set $A_r=A\cap(r+b\Z)$.
    Additionally set $P_r=(A_r\pm b)\setminus A$.
    Note that $P_r\subseteq r+b\Z$ and so the $P_r$'s are pairwise disjoint and, certainly, $P_r\subseteq(A\pm B)\setminus A$.
    Since $A_r\neq\varnothing$ and is finite, it is easy to see that $\abs{P_r}\geq 2$ with equality if and only if $A_r$ is an arithmetic progression with step size $b$.
    In particular,
    \begin{equation}\label{eqn:boundaryP}
        \abs{(A\pm B)\setminus A}\geq\sum_{r=0}^{b-1}\abs{P_r}\geq 2b.
    \end{equation}

    To finish the claim, suppose that equality holds in \cref{eqn:boundaryP} but that $A$ is not an interval.
    In order for the latter inequality of \cref{eqn:boundaryP} to hold with equality, it must be the case that $\abs{P_r}=2$, and thus $A_r$ is an arithmetic progression with step size $b$, for each $r\in[0,b-1]$.

    For the sake of brevity, in what follows, all subscripts are computed modulo $b$; that is, $A_x=A_{y}$ and $P_x = P_y$ for any $x \equiv y \pmod b$.
    \medskip

    Let $x$ be the smallest integer larger than $\min A$ which does not belong to $A$, so $x-1\in A$.
    Thus, $x=(x-1)+1\in(A\pm B)\setminus A$.
    If $x\notin P_x$, then $\abs{(A\pm B)\setminus A}\geq 2b+1$ which yields our contradiction, so we may suppose that $x\in P_x$.
    Since $A_x$ is an arithmetic progression with step size $b$, either $x=\min A_x-b$ or $x=\max A_x+b$.

    Symmetrically, let $y$ be the largest integer smaller than $\max A$ which does not belong to $A$, so $y+1\in A$.
    By symmetric logic to that in the previous paragraph, we must have either $y=\min A_y-b$ or $y=\max A_y+b$ or else we have reached a contradiction.
    Furthermore, since $A$ is not an interval, observe that $x\leq y$.
    From here we break into cases.
    \medskip

    \textbf{Case 1:} $x=\min A_x-b$ or $y=\max A_y+b$.
    These two situations are symmetric and so we consider only the case when $x=\min A_x-b$.

    Set $z=\min A_{x-1}$; since $x-1\in A$ we must have $z\leq x-1$.
    Now, since $z-(b-1)<x$ we know that $z-(b-1)$ is not a member of any $P_r$ since $z-(b-1)\equiv x\pmod b$.
    In addition to the previous observation we know that $[\min A,x-1]\subseteq A$, and so $z-(b-1)\notin A$.
    This implies that $z-(b-1)\in(A\pm B)\setminus A$, and so $\abs{(A\pm B)\setminus A}\geq 2b+1$, yielding a contradiction.
    \medskip

    \textbf{Case 2:} $x=\max A_x+b$ and $y=\min A_y-b$.
    By the definition of $x$ and $y$, we know that $[\min A,x-1]\cup[y+1,\max A]\subseteq A$.
    Thus, in this case we have $x\geq \min A+b$ and $y\leq\max A-b$.
    In particular, both $[\min A,x-1]$ and $[y+1,\max A]$ contain every residue modulo $b$.
    However, $y\geq x>\max A_x$ and so $[y+1,\max A]\subseteq A\setminus A_x$ must be void of any integers which are $x$ modulo $b$; a contradiction.
\end{proof}

We next strengthen the previous observation, showing that the same conclusion holds if only a single residue class is absent.

\begin{lemma}\label[lemma]{missing_one}
    Fix a finite set $B\subseteq\Z_{\geq 1}$, set $b=\max B$ and suppose that $1,b-1\in B$ as well.
    If $A\subseteq\Z$ is a finite set which contains all but one residue modulo $b$, then $\abs{(A\pm B)\setminus A}\geq 2b$ with equality if and only if $A$ is an interval.
\end{lemma}
\begin{proof}
    For each $r\in[0,b-1]$, set $A_r=A\cap(r+b\Z)$.
    Additionally set $P_r=(A_r\pm b)\setminus A$.
    By translating if necessary, we may suppose that $A_0=\varnothing$ while $A_r\neq\varnothing$ for all other values of $r$.
    Again, the $P_r$'s are pairwise disjoint and $P_r\subseteq(A\pm B)\setminus A$.
    Additionally, $\abs{P_r}\geq 2$ for each $r\in[b-1]$.

    Now, additionally consider $Q_1=A_1+\{-1,b-1\}$ and $Q_{b-1}=A_{b-1}+\{1,-(b-1)\}$.
    Observe that $Q_1,Q_{b-1}\subseteq b\Z$ and so $Q_1,Q_{b-1}$ are disjoint from each of $P_1,\dots,P_{b-1}$; however, they may not be disjoint from each other.
    Furthermore, $Q_1,Q_{b-1}\subseteq(A\pm B)\setminus A$ since $A_0=\varnothing$.
    Although $Q_1,Q_{b-1}$ may not be disjoint from one another, $\abs{Q_1\cup Q_{b-1}}\geq 2$ since each of these sets has size at least $2$.
    Thus, all together, we have
    \[
        \abs{(A\pm B)\setminus A}\geq\sum_{r=1}^{b-1}\abs{P_r}+\abs{Q_1\cup Q_b{-1}}\geq 2(b-1)+2=2b.
    \]
    If equality holds, then certainly $\abs{Q_1\cup Q_{b-1}}=2$ and so $Q_1=Q_{b-1}$.
    This implies that $A_1=\{a_1\}$ and $A_{b-1}=\{a_{b-1}\}$ where $a_1+(b-1)=a_{b-1}+1$.
    Translating $A$ by some multiple of $b$ if necessary, we may assume that $a_1=1$, so then $a_{b-1}=b-1$.
    Now, since $A$ contains every residue modulo $b$ except for $0$, the only way for $A$ to not be an interval is if $A$ contains an element outside of $[b-1]$.
    Suppose that $\max A>b-1$ (the situation where $\min A<1$ follows symmetrically); note that actually $\max A>b+1$ since $b,b+1\notin A$.
    Let $x\in[b,\max A]$ be the largest integer which is \emph{not} a member of $A$.
    Then $x+1\in A$ and so $x+b=(x+1)+(b-1)\in A\pm B$.
    Additionally, $x+b>\max A$ since if $x+b\leq\max A$, then $[x+1,x+b]\subseteq A$ which is impossible since this interval contains a multiple of $b$.
    Finally, $x+b$ does not belong to any $P_r$ since $x\notin A$ and also $x+b\notin Q_1\cup Q_{b-1}$ since $x\geq b+1$.
    Thus, $\abs{(A\pm B)\setminus A}\geq 2b+1$ if $A$ is not an interval.
\end{proof}

The previous two lemmas come together to yield the statement we are actually after.

\begin{theorem}\label[theorem]{two_consecutive}
    Fix $b\geq 3$, set $B=\{1,b-1,b,b+1\}$ and let $A\subseteq\Z$ be some finite set.
    If $A$ is \emph{not} missing two consecutive residues modulo $b$, then $\abs{(A\pm B)\setminus A}\geq 2(b+1)$ with equality if and only if $A$ is an interval.
\end{theorem}
\begin{proof}
    Set $B'=\{1,b-1,b\}$.
    Certainly $\min A-(b+1)$ and $\max A+(b+1)$ belong to $(A\pm B)\setminus A$, but neither belongs to $(A\pm B')\setminus A$.
    Consequently, $\abs{(A\pm B)\setminus A}\geq 2+\abs{(A\pm B')\setminus A}$.
    As such, if $A$ is missing at most one residue modulo $b$, then the claim follows from \Cref{missing_one,missing_none}.

    Therefore, in order to conclude the proof, we must show that if $A$ is missing at least two residues modulo $b$, yet no two consecutive residues, then $\abs{(A\pm B)\setminus A}>2(b+1)$ since it is impossible for $A$ to be an interval under these circumstances.
    \medskip

    Let $R\subseteq[0,b-1]$ be the set of residues modulo $b$ contained in $A$.
    By translating if necessary, we may suppose that $0\in R$ and $b-1\notin R$, so $R\subseteq[0,b-2]$.
    Observe that we can write $R=I_1\cup\dots\cup I_k$ where $I_1,\dots,I_k$ are disjoint intervals and $\max I_i+1\notin R$ for each $i\in[k]$.
    Since two residues are missing, but no two consecutive residues are, observe that $k\geq 2$ and that $k+\abs R=b$.
    \medskip

    Similarly to the proofs of \Cref{missing_one,missing_none}, for each $r\in R$, define $A_r=A\cap(r+b\Z)$ and set $P_r=(A_r\pm b)\setminus A$.
    As before, $\abs{P_r}\geq 2$ for each $r\in R$ and they are pairwise disjoint.
    Now, for each $i\in[k]$, define $Q_i=A_{m_i}+\{1,b+1,-(b-1)\}$ where $m_i=\max I_i$.
    Note that $Q_i\subseteq (m_i+1)+b\Z$ and so $Q_1,\dots,Q_k$ are pairwise disjoint and also disjoint from $A$ and $P_r$ for all $r\in R$ since $m_i+1\notin R$.
    Furthermore, $\abs{Q_i}\geq 3$ for each $i\in[k]$ since $A_{m_i}$ is non-empty.
    Finally, the element $\min A-(b+1)$ does not belong to any $P_r$ nor $Q_i$ and so
    \[
        \abs{(A\pm B)\setminus A}\geq \sum_{r\in R}\abs{P_r}+\sum_{i=1}^k\abs{Q_i}+1\geq 2\abs R+3k+1=2b+k+1\geq 2b+3.\qedhere
    \]
\end{proof}

We can now perform very similar reasoning to \Cref{edge_not_nested} in order to establish a phase-transition in the vertex-isoperimetric problem on $\cay(\Z,\pm\{1,b-1,b,b+1\})$.
Recalling the discussion at the beginning of this section, the key idea is to take an optimal $A$, which is not an interval, and ``cut'' along some consecutive missing residue classes modulo $b$ in order to embed $A$ into $\ell_\infty^2$ without affecting its boundary.

\begin{theorem}\label[theorem]{vert_not_nested}
    Fix a positive integer $b\geq 3$ and set $G=\cay(\Z,\pm\{1,b-1,b,b+1\})$.
    For every integer $n\geq b$:
    \begin{enumerate}
        \item If $\sqrt n>{b-1\over 2}$, then every element of $\Opt{v}{G}(n)$ is an interval of length $n$.
        \item If $n=k^2$ and $k<{b-1\over 2}$, then every element of $\Opt{v}{G}(n)$ is a translate of $[k]+b[k]$.
    \end{enumerate}
\end{theorem}
\begin{proof}
    Since $n\geq b$, observe that
    \begin{equation}\label{eqn:vert_interval}
        \abs{\boundary{v}{G}[n]}=2(b+1).
    \end{equation}

    Suppose that $A\in\Opt{v}{G}(n)$ is \emph{not} an interval.
    We will show the following:
    \begin{enumerate}
        \item $\sqrt n\leq{b-1\over 2}$, and
        \item If $n=k^2$ and $k<{b-1\over 2}$, then $\abs{\boundary{v}{G}A}<\abs{\boundary{v}{G}[n]}$ and $A$ is a translate of $[k]+b[k]$.
    \end{enumerate}
    This will establish the claim.
    \medskip

    Thanks to \Cref{two_consecutive}, since $A$ is not an interval, we know that $A$ must be missing two consecutive residues modulo $b$.
    Without loss of generality, we may suppose that these two residues are $0$ and $b-1$.
    \medskip

    Consider the map $\phi\colon\Z\to\Z^2$ defined by $\phi(x)=(q,r)$ where $x=qb+r$ is the Euclidean decomposition of $x$; note that the image of $\phi$ is $\Z\times[0,b-1]$.
    Since $A$ has no elements which are $0$ or $b-1$ modulo $b$, we know that $\phi(A)\subseteq\Z\times[1,b-2]$.
    Because of this, for any $a\in A$, we have
    \begin{align*}
        \phi(a\pm 1) &=\phi(a)\pm \e_2,\\
        \phi(a\pm(b-1)) &=\phi(a)\pm(\e_1-\e_2),\\
        \phi(a\pm b) &=\phi(a)\pm \e_1,\\
        \phi(a\pm(b+1)) &=\phi(a)\pm(\e_1+\e_2).
    \end{align*}
    In particular, $\phi\bigl(\boundary{v}{G}A\bigr)=\boundary{v}{\ell_\infty^2}\phi(A)$.
    As a consequence of the vertex-isoperimetric inequality for $\ell_\infty^2$ (\Cref{infty_iso}), we have $\abs{\boundary{v}{G}A}\geq 4(\sqrt n+1)$ with equality if and only if $\phi(A)$ is a translate of $[k]^2$.

    From \cref{eqn:vert_interval}, this implies that $4(\sqrt n+1)\leq 2(b+1)\implies \sqrt n\leq{b-1\over 2}$, which establishes the first part of the claim.

    Now, suppose that $n=k^2$ and that $k=\sqrt n<{b-1\over 2}$.
    In particular, $k<b-1$ and so $[k]^2\subseteq\Z\times[1,b-2]$.
    Thus, from the observation above, we must have $\abs{\boundary{v}{G}A}=4(k+1)$ with equality if and only if $\phi(A)$ is a translate of $[k]^2$.
    Of course, $\phi(A)$ is a translate of $[k]^2$ if and only if $A$ is a translate of $[k]+b[k]$.
    This concludes the proof since here $4(k+1)<2(b+1)=\abs{\boundary{v}{G}[n]}$.
\end{proof}

\subsection{Proof of \texorpdfstring{\Cref{not_nested}}{Theorem 1.2}}\label{sec:not_nested}

Given a positive integer $N$, fix any even integer $b\geq \max\{N,10\}$ and consider the Cayley graphs $G_{\operatorname e}=\cay(\Z,\pm\{1,b\})$ and $G_{\operatorname v}=\cay(\Z,\pm\{1,b-1,b,b+1\})$.

Set $k={b-2\over 2}$, $n_1=k^2$ and $n_2=(b+1)(k-1)$.
The following can be checked since $b\geq\max\{N,10\}$: $n_1\geq N$, and $n_2-n_1\geq N$, and $k=\sqrt{n_1}<{b-1\over 2}$ and $\sqrt{n_2}>{b+1\over 2}$.

Now, thanks to \Cref{edge_not_nested,vert_not_nested}, we know that $\Opt{e}{G_e}(n_1)=\Opt{v}{G_v}(n_1)$ is the set of all translates of $[k]+b[k]$.
Additionally, $\Opt{e}{G_e}(n_2)=\Opt{v}{G_v}(n_2)$ is the set of all intervals of length $n_2$.
Now, if $A$ is any translate of $[k]+b[k]$, then $\max A-\min A=(b+1)(k-1)>n_2-1$; thus $A$ is not contained in an interval of length $n_2$, which concludes the proof.
\qed

\section{Intervals are eventually optimal}\label{sec:interval_opt}

We conclude this paper by showing that intervals are eventually optimal for both the vertex- and the edge-isoperimetric problems on any Cayley graph on $\Z$.
\medskip

Recall that $B\subseteq\Z$ generates $\Z$ if and only if $\cay(\Z,B)$ is weakly connected.
We begin with a straightforward structural result about such Cayley graphs.

Say that graph $G$ is \emph{rooted} a vertex $r$ if there is a directed path from $r$ to $v$ for every other vertex $v\in V(G)$.
Note that $G$ is strongly-connected if and only if it is rooted at every vertex.

\begin{prop}\label[prop]{rootedResidues}
    Fix a finite set $B\subseteq\Z$.
    The following are equivalent:
    \begin{enumerate}
        \item $B$ generates $\Z$.
        \item For every $n\in\Z_{\geq 1}$, there is a finite subgraph of $\cay(\Z,B)$ which is rooted at $0$ and contains an element of every residue class modulo $n$.
    \end{enumerate}
\end{prop}
\begin{proof}
    (2)$\implies$(1): Fix any $x\in\Z$.
    Without loss of generality, $B$ contains a positive integer: call one of them $b_1$.
    By assumption, there is a directed path in $\cay(\Z,B)$ from $0$ to some element of $x+b_1\Z$.
    In particular, if $B=\{b_1,\dots,b_k\}$ then we can write $qb_1+x=n_1b_1+\dots+n_kb_k$ for some $q\in\Z$ and some $n_1,\dots,n_k\in\Z_{\geq 0}$.
    Therefore, $x=(n_1-q)b_1+n_2b_2+\dots+n_kb_k$ and so $x$ can be written as an integer combination of elements of $B$.
    We conclude that $B$ generates $\Z$.
    \medskip

    (1)$\implies$(2): Write $B=\{b_1,\dots,b_k\}$ and fix any $x\in\Z$.
    Since $B$ generates $\Z$, we can write $x=n_1b_1+\dots+n_kb_k$ where $n_1,\dots,n_k\in\Z$.
    Set $I=\{i\in[k]:n_i<0\}$ and set $N=-\sum_{i\in I}n_ib_i$.
    Thus, for any $n\in \Z_{\geq 1}$, we have
    \begin{align*}
        x+nN &=(n_1b_1+\dots+n_kb_k)-\sum_{i\in I}n\cdot n_ib_i =\sum_{i\in I}(1-n)n_ib_i+\sum_{i\notin I}n_ib_i.
    \end{align*}
    Since $n\geq 1$, this expresses $x+nN$ as a non-negative-integer combination of the elements of $B$.
    In other words, there is a directed path from $0$ to $x+nN$ in $\cay(\Z,B)$.
    Calling this path $P_x$, we find that the graph $P_0\cup P_1\cup\dots\cup P_{n-1}$ is a finite subgraph of $\cay(\Z,B)$ rooted at $0$ which contains an element of every residue class modulo $n$.
\end{proof}

Using this structure, we next show that any set $A\subseteq\Z$ which is void of some residue modulo an element of $B$ must have large vertex-boundary in $\cay(\Z,B)$.

\begin{lemma}\label[lemma]{missing_implies_large}
    Fix a finite set $B$ which generates $\Z$.
    There is some $\epsilon=\epsilon(B)>0$ such that, for every finite $A\subseteq\Z$, either:
    \begin{enumerate}
        \item $A$ contains every residue modulo $\abs b$ for each $b\in B$, or
        \item $\abs{(A+B)\setminus A}\geq\epsilon\abs A$.
    \end{enumerate}
\end{lemma}
\begin{proof}
    Suppose that $A$ is void of some residue modulo $\abs b$ for some $b\in B$.
    Of course, we may suppose that $\abs b\geq 2$ since otherwise $A$ must be empty.
    Since $B$ generates $\Z$, let $T$ be any finite subgraph of $\cay(\Z,B)$ which is rooted at $0$ and contains an element of every residue class modulo $\abs b$.
    The existence of such a $T$ is guaranteed by \Cref{rootedResidues}.
    We treat $T$ as both a graph and as a simple subset of $\Z$ in what follows.

    Fix any $a\in A$ and consider the graph $a+T$.
    Of course, the graph $a+T$ is rooted at $a$ and still contains an element of every residue class modulo $\abs b$.
    Thus, since $A$ is void of some residue modulo $\abs b$, not every element of $a+T$ belongs to $A$.
    If, say, $z\in (a+T)\setminus A$, then the directed path from $a$ to $z$ in $a+T$ must contain an element of $(A+B)\setminus A$ by construction.
    Therefore, for every $a\in A$, there is some $z\in(A+B)\setminus A$ with $z\in a+T$.
    We conclude that
    \[
        \abs{(A+B)\setminus A}\geq{\abs A\over\max T-\min T},
    \]
    since $B$ is finite and $T$ depends only on $b$ and $B$.
\end{proof}

We can now prove that intervals are always eventually optimal for both vertex- and edge-isoperimetric problems on Cayley graphs on $\Z$.
\begin{proof}[Proof of \Cref{interval_opt}]
    Fix a finite set $B$ which generates $\Z$, set $G=\cay(\Z,B)$ and let $\epsilon=\epsilon(B)>0$ be as in \Cref{missing_implies_large}.

    We wish to show that for all $n$ sufficiently large, $\abs{\boundary{e}{G}A}\geq\abs{\boundary{e}{G}[n]}$ and $\abs{\boundary{v}{G}A}\geq\abs{\boundary{v}{G}[n]}$ whenever $A\subseteq\Z$ has size $n$.

    \paragraph{The edge-isoperimetric inequality on $G$.}
    Consider a large integer $n$ and suppose that $A\subseteq\Z$ is any set of size $n$.
    Observe that if $n\geq\max_{b\in B}\abs b$, then $\abs{\boundary{e}{G}[n]}=\sum_{b\in B}\abs b$ (\Cref{every_res_then_interval}).
    \medskip

    \textbf{Case 1:} $A$ contains every residue modulo $\abs b$ for each $b\in B$.
    Then, due to \Cref{every_res_then_interval}, $\abs{\boundary{e}{G}A}\geq\abs{\boundary{e}{G}[n]}$, as needed.
    \medskip

    \textbf{Case 2:} $A$ is void of some residue modulo $\abs b$ for some $b\in B$.
    Then $\abs{\boundary{e}{G}A}\geq\abs{(A+B)\setminus A}\geq\epsilon n$ due to \Cref{missing_implies_large}.
    Of course, if $n$ is sufficiently large, then $\epsilon n>\sum_{b\in B}\abs b=\abs{\boundary{e}{G}[n]}$ as needed.

    \paragraph{The vertex-isoperimetric inequality on $G$.}
    Begin by decomposing $B=B^+\cup(-B^-)$ where $B^+,B^-\subseteq\Z_{\geq 1}$.
    Set $b^+=\max B^+$ and $b^-=\max B^-$ where here we define $\max\varnothing=0$.
    Observe that if $n\geq \max\{b^+,b^-\}$, then $\abs{\boundary{v}{G}[n]}=b^++b^-$.
    Without loss of generality, we may suppose that $B^+\neq\varnothing$ since we may replace $B$ by $-B$ without changing the problem.
    \medskip

    Consider a large integer $n$ and suppose that $A\subseteq\Z$ is any set of size $n$.
    \medskip

    \textbf{Case 1:} $A$ is void of some residue modulo $\abs b$ for some $b\in B$.
    Then $\abs{\boundary{v}{G}A}=\abs{(A+B)\setminus A}\geq\epsilon n$ due to \Cref{missing_implies_large}.
    Of course, if $n$ is sufficiently large, then $\epsilon n>b^++b^-=\abs{\boundary{v}{G}[n]}$ as needed.
    \medskip

    \textbf{Case 2:} $A$ contains every residue modulo $\abs b$ for each $b\in B$.
    \medskip

    \textbf{Case 2a:} $B^-=\varnothing$, and so $b^-=0$.
    Here we see that $\bigl\{\max\bigl(A\cap(r+b^+\Z)\bigr)+b^+:r\in[0,b^+-1]\bigr\}$ is a set of $b^+$ distinct elements of $\boundary{v}{G}A$ and so $\abs{\boundary{v}{G}A}\geq b^+=\abs{\boundary{v}{G}[n]}$.
    \medskip

    \textbf{Case 2b:} $B^-\neq\varnothing$.
    For each $r\in[0,b^+-1]$, define $a^+_r=\max\bigl(A\cap(r+b^+\Z)\bigr)$.
    For each $s\in[0,b^--1]$, define $a^-_s=\min\bigl(A\cap(s+b^-\Z)\bigr)$.
    Of course, $a^+_0+b^+,a^+_1+b^+,\dots,a^+_{b^+-1}+b^+$ are distinct elements of $(A+B)\setminus A$ and similarly $a^-_0-b^-,a^-_1-b^-,\dots,a^-_{b^--1}-b^-$ are also distinct elements of $(A+B)\setminus A$.
    If all of these $b^++b^-$ elements are distinct, then $\abs{\boundary{v}{G}A}\geq b^++b^-=\abs{\boundary{v}{G}[n]}$ and we have established the claim.

    Thus, suppose that $a^+_r+b^+=a^-_s-b^-$ for some $r\in[0,b^+-1]$ and $s\in[0,b^--1]$.
    Of course, this would mean that $A^+=A\cap[a^-_s,\infty)$ has no element which is $r$ modulo $b^+$ and that $A^-=A\cap(-\infty,a^+_r]$ has no element which is $s$ modulo $b^-$.
    Therefore, $\abs{\boundary{v}{G}A^+}\geq\epsilon\abs{A^+}$ and $\abs{\boundary{v}{G}A^-}\geq\epsilon\abs{A^-}$ due to \Cref{missing_implies_large}.

    Since $a^-_s-a^+_r=b^++b^-$, we see that $\abs{A^+}+\abs{A^-}\geq\abs A-(b^++b^-)$.
    Finally, it is clear that $\boundary{v}{G}A^+\cap\boundary{v}{G}A^-=\{a^+_r+b^+=a^-_s-b^-\}$ and that $\boundary{v}{G}A^+\cup\boundary{v}{G}A^-\subseteq\boundary{v}{G}A\cup[a^+_r,a^-_s]$.
    We conclude that
    \begin{align*}
        \abs{\boundary{v}{G}A} &\geq\abs{\boundary{v}{G}(A^+)}+\abs{\boundary{v}{G}(A^-)}-b^+-b^--1\\
                                &\geq\epsilon\abs{A^+}+\epsilon\abs{A^-}-b^+-b^--1\\
                                &\geq \epsilon n- 2(b^++b^-+1).
    \end{align*}
    Of course, for $n$ sufficiently large, $\epsilon n-2(b^++b^-+1)\geq b^++b^-=\abs{\boundary{v}{G}[n]}$, which establishes the claim.
\end{proof}

\section{Concluding remarks}\label{sec:remarks}

While we have built counterexamples to the original question of Barber--Erde, our constructions do not conflict with the heart of their question.
In particular, our counterexamples disappear if one is allowed to ignore finitely many optimal sets.
It is therefore natural to ask:
\begin{question}\label[question]{q.eventuallynested}
        Let $B\subseteq\Z^d$ be a finite set which generates $\Z^d$ and set $G=\cay(\Z^d,B)$.
    \begin{itemize}
        \item Does there exist an integer $N$ and an ordering $u_1,u_2,\ldots$ of $\Z^d$ for which $\{u_1,\dots,u_n\}\in\Opt{v}{G}(n)$ for every $n\geq N$?
    \end{itemize}
\end{question}
\Cref{interval_opt} answers the above question affirmatively when $d=1$. Note that Question 2 in \cite{barber_iso} asks whether their asymptotically-best constructions are also in $\Opt{\mathnormal i}{G}(n)$ for infinitely many values of $n$. While this is not an immediate strengthening or weakening of the above question, it is likely that any method that resolves one will also resolve the other.

We also remark that, in an earlier version of this paper, we also asked the corresponding edge-isoperimetric version of \Cref{q.eventuallynested}. In fact, this has since already been answered \emph{in the negative} by Strachan and Swanepoel \cite{strachan2025edge}.

\medskip

Our second question is much less general and is more of a curiosity related to our constructions.

Consider the toroidal-cylindrical grids
\begin{align*}
    C_1^k(b) &\eqdef\cay\bigl(\Z\times(\Z/b\Z)^k,\pm\{\e_1,\dots,\e_{k+1}\}\bigr),\qquad\text{and}\\
    C_\infty^k(b) &\eqdef\cay\bigl(\Z\times(\Z/b\Z)^k,\{0,\pm 1\}^{k+1}\setminus\{\mathbf 0\}\bigr).
\end{align*}
While we are not aware of any explicit results in this direction, the edge-isoperimetric problem for $C_1^k(b)$ and the vertex-isoperimetric problems for $C_\infty^k(b)$ should each exhibit $k$ distinct phase transitions---it seems likely that the techniques from \cite{bollobas_edge} can be used to prove this.
Following the motivation in \Cref{sec:motivation}, consider the map $\phi\colon\Z\to\Z\times(\Z/b\Z)^k$ defined by $\phi(x)=(x_k,x_{k-1},\dots,x_0)$ where $x=x_k b^k+x_{k-1}b^{k-1}+\dots+x_0$ with $x_0,\dots,x_{k-1}\in[0,b-1]$.
This map satisfies $\phi(b^i)=\e_{k+1-i}$ for each $i\in[0,k]$.
It is natural to wonder if the constructions in this paper can be extended:
\begin{question}
    Fix positive integers $b,k$.
    If $b$ is sufficiently large compared to $k$:
    \begin{itemize}
        \item Does the edge-isoperimetric problem on $\cay\bigl(\Z,\pm\{1,b,b^2,\dots,b^k\}\bigr)$ exhibit $k$ distinct phase transitions?
        \item Does the vertex-isoperimetric problem on $\cay\bigl(\Z,\bigl\{\epsilon_0+\epsilon_1 b+\dots+\epsilon_k b^k:\epsilon_0,\dots,\epsilon_k\in\{0,\pm 1\}\bigr\}\bigr)$ exhibit $k$ distinct phase transitions?
    \end{itemize}

\end{question}

\appendix

\section{The vertex-isoperimetric inequality on \texorpdfstring{$\ell_\infty^2$}{the infinity-grid}}\label{sec:isoproofs}
Here we prove \Cref{infty_iso}.

For a set $A\subseteq\Z^2$ and integer $i\in\Z$, define $A_i\eqdef\{x\in\Z:(x,i)\in A\}$.
Observe that $A=\bigcup_{i\in\Z}\bigl(A_i\times\{i\}\bigr)$.
Additionally define $N[A]\eqdef A+\{0,\pm 1\}^2=A\cup\boundary{v}{\ell_\infty^2}A$, so $\abs{N[A]}=\abs A+\abs{\boundary{v}{\ell_\infty^2}A}$.
It is not difficult to see that $N[A]_i=(A_{i-1}\cup A_i\cup A_{i+1})+\{0,\pm 1\}$ for each $i\in\Z$.
Finally, define
\[
    C(A)\eqdef\bigcup_{i\in\Z}\bigl([a_i]\times\{i\}\bigr),
\]
where $a_i=\abs{A_i}$.
Of course, $C$ is idempotent and $\abs{C(A)}=\abs A$ for any $A\subseteq\Z^2$.
\begin{lemma}\label[lemma]{nestedGetsSmaller}
    Let $A\subseteq\Z^2$ be a finite set.
    Then
    \begin{enumerate}
        \item $\abs{N[C(A)]}\leq\abs{N[A]}$, and
        \item If $C(A)$ is a translate of $[k]^2$ for some integer $k$, then $\abs{N[C(A)]}=\abs{N[A]}$ if and only if $A$ is also a translate of $[k]^2$.
    \end{enumerate}
\end{lemma}
\begin{proof}
    For each $i\in\Z$, set $a_i=\abs{A_i}$.
    Define $\mathcal I=\{i\in\Z:N[A]_i\neq\varnothing\}$.
    Certainly if $X\subseteq\Z$ is finite and non-empty, then $\abs{X+\{0,\pm 1\}}\geq 2+\abs X$ with equality if and only if $X$ is an interval.
    Therefore,
    \begin{equation}\label{infinitySlices}
        \abs{N[A]}=\sum_{i\in\mathcal I}\abs{N[A]_i}=\sum_{i\in\mathcal I}\abss{(A_{i-1}\cup A_i\cup A_{i+1})+\{0,\pm 1\}}\geq\sum_{i\in\mathcal I}\bigl(2+\abs{A_{i-1}\cup A_i\cup A_{i+1}}\bigr),
    \end{equation}
    with equality if and only if $A_{i-1}\cup A_i\cup A_{i+1}$ is an interval for all $i\in\mathcal I$.
    On the other hand, since $C(A)_i=[a_i]$ for all $i\in\Z$,
    \begin{align*}
        \abs{N[C(A)]} &=\sum_{i\in\mathcal I}\abs{N[C(A)]_i}=\sum_{i\in\mathcal I}\bigl(2+\abs{C(A)_{i-1}\cup C(A)_i\cup C(A)_{i+1}}\bigr)\\
                      &=\sum_{i\in\mathcal I}\bigl(2+\max\{a_{i-1},a_i,a_{i+1}\}\bigr).
    \end{align*}
    Of course, $\abs{A_{i-1}\cup A_i\cup A_{i+1}}\geq\max\{a_{i-1},a_i,a_{i+1}\}$ and so we have shown that $\abs{N[C(A)]}\leq\abs{N[A]}$, thus establishing the first part of the claim.

    To prove the second part of the claim, we note that $\abs{N[C(A)]}=\abs{N[A]}$ if and only if for each $i\in\mathcal I$, we have that $A_{i-1}\cup A_i\cup A_{i+1}$ is an interval and that $\abs{A_{i-1}\cup A_i\cup A_{i+1}}=\max\{a_{i-1},a_i,a_{i+1}\}$.
    Of course, if $C(A)$ is a translate of $[k]^2$, then we must have $C(A)=[k] \times [I,J]$ where $J-I+1=k$.
    In particular, $\mathcal I=[I-1,J+1]$ and $a_i=k$ for each $i\in[I,J]$.
    Now, considering $i=I-1$, we find that $A_{i+1}=A_I$ must be an interval of length $k$.
    We claim that, in fact, $A_i=A_I$ for all $i\in[I,J]$, which will establish the claim.
    Indeed, suppose for the sake of contradiction that $I<i\leq J$ is the smallest integer for which $A_i\neq A_I$.
    Since the three sets $A_{i-1}$, $A_i$ and $A_{i-1}\cup A_i\cup A_{i+1}$ each have size $k$, we must have $A_i=A_{i-1}=A_I$ as needed.
\end{proof}

Say that a set $A\subseteq\Z^2$ is \emph{nested} if there are integers $I\leq J$ for which
\begin{itemize}
    \item $A_i=\varnothing$ for all $i\notin[I,J]$, and
    \item $A_I\supseteq \dots\supseteq A_J\neq\varnothing$.
\end{itemize}

\begin{lemma}\label[lemma]{nestedAreSquares}
    If $A\subseteq\Z^2$ is finite and nested, then $\abs{\boundary{v}{\ell_\infty^2}A}\geq 4\bigl(\sqrt{\abs A}+1\bigr)$ with equality if and only if $A$ is a translate of $[k]^2$ for some positive integer $k$.
\end{lemma}
\begin{proof}
    Using \cref{infinitySlices} and the fact that $A$ is nested, we bound
    \begin{align*}
        \abs{N[A]} &\geq\sum_{i=I-1}^{J+1}\bigl(2+\abs{A_{i-1}\cup A_i\cup A_{i+1}}\bigr)=\bigl(2+\abs{A_I}\bigr)+\bigl(2+\abs{A_I}\bigr)+\sum_{i=I+1}^{J+1}\bigl(2+\abs{A_{i-1}}\bigr)\\
                   &=4+2(J-I+1)+2\abs{A_I}+\sum_{i=I}^J\abs{A_i}=4+2(J-I+1)+2\abs{A_I}+\abs A\\
                   &\geq 4\bigl(\sqrt{(J-I+1)\abs{A_I}}+1\bigr)+\abs A\geq 4\bigl(\sqrt{\abs A}+1\bigr)+\abs A.
    \end{align*}
    The first inequality is tight if and only if each $A_i$ is an interval.
    The second inequality follows from the arithmetic--geometric mean inequality and is tight if and only if $J-I+1=\abs{A_I}$.
    The last inequality follows from the fact that $\abs A=\sum_{i=I}^J\abs{A_I}\leq (J-I+1)\abs{A_I}$, and is tight if and only if $\abs{A_I}=\dots=\abs{A_J}$.
    Note that $\abs{A_I}=\dots=\abs{A_J}$ if and only if $A_I=\dots=A_J$ which means that $A=[I,J]\times A_I$.

    Thus, we have shown that $\abs{N[A]}\geq 4\bigl(\sqrt{\abs A}+1\bigr)+\abs A$ with equality if and only if $A$ is a translate of $[k]^2$ for some positive integer $k$.
    Since $\abs{N[A]}=\abs A+\abs{\boundary{v}{\ell_\infty^2}A}$, the claim follows.
\end{proof}

We can finally prove \Cref{infty_iso}.

For a set $A\subseteq\Z^2$, define the transpose of $A$ to be $A^T\eqdef\{(x,y)\in\Z^2:(y,x)\in A\}$.
Of course, $\boundary{v}{\ell_\infty^2}A^T=\bigl(\boundary{v}{\ell_\infty^2}A\bigr)^T$ and so transposing a set does not affect the size of its boundary.
\begin{proof}[Proof of \Cref{infty_iso}]
    Fix any set $A\subseteq\Z^2$ of size $n$.
    Then, thanks to the first part of \Cref{nestedGetsSmaller}, we have
    \begin{equation}\label{compressionIneq}
        \abs{\boundary{v}{\ell_\infty^2}A}\geq\abs{\boundary{v}{\ell_\infty^2}C(A)}=\abs{\boundary{v}{\ell_\infty^2}C(A)^T}\geq\abs{\boundary{v}{\ell_\infty^2}C\bigl(C(A)^T\bigr)}.
    \end{equation}
    It can be checked that $C\bigl(C(A)^T\bigr)$ is nested and so $\abs{\boundary{v}{\ell_\infty^2}C\bigl(C(A)^T\bigr)}\geq 4(\sqrt n+1)$ with equality if and only if $C\bigl(C(A)^T\bigr)$ is a translate of $[k]^2$ for some positive integer $k$ due to \Cref{nestedAreSquares}.
    Therefore, $\abs{\boundary{v}{\ell_\infty^2}A}\geq 4(\sqrt n+1)$.

    Now, if $\abs{\boundary{v}{\ell_\infty^2}A}=4(\sqrt n+1)$, then every inequality in \cref{compressionIneq} holds with equality and $C\bigl(C(A)^T\bigr)$ is a translate of $[k]^2$.
    Then, thanks to the second part of \Cref{nestedGetsSmaller}, we must have that $C(A)^T$ (and thus $C(A)$) is also a translate of $[k]^2$.
    Applying the second part of \Cref{nestedGetsSmaller} once more allows us to conclude that $A$ must also be a translate of $[k]^2$, which establishes the claim.
\end{proof}

\bibliographystyle{abbrv}
\bibliography{references}

@article{barber_iso,
	journal={Discrete Analysis},
	doi={10.19086/da.3555},
	title={Isoperimetry in integer lattices},
	author={Barber, Ben and Erde, Joshua},
	date={2018-04-20},
	year=2018,
	month=4,
	day=20,
}

@article{barber_stability,
  title={Isoperimetric stability in lattices},
  author={Barber, Ben and Erde, Joshua and Keevash, Peter and Roberts, Alexander},
  journal={Proceedings of the American Mathematical Society},
  volume={151},
  number={12},
  pages={5021--5029},
  year={2023}
}

@article{bollobas_edge,
  doi = {10.1007/bf01275667},
  url = {https://doi.org/10.1007/bf01275667},
  year = {1991},
  month = dec,
  publisher = {Springer Science and Business Media {LLC}},
  volume = {11},
  number = {4},
  pages = {299--314},
  author = {B{\'e}la Bollob{\'a}s and Imre Leader},
  title = {Edge-isoperimetric inequalities in the grid},
  journal = {Combinatorica}
}

@article{gupta_planar,
  title={Planar Lattice Subsets with Minimal Vertex Boundary},
  author={Gupta, Radhika and Levcovitz, Ivan and Margolis, Alexander and Stark, Emily},
  journal={The Electronic Journal of Combinatorics},
  volume={28},
  number={3},
  year={2021}
}

@article{lev_edge,
  title={Edge-isoperimetric problem for Cayley graphs and generalized Takagi functions},
  author={Lev, Vsevolod F},
  journal={SIAM Journal on Discrete Mathematics},
  volume={29},
  number={4},
  pages={2389--2411},
  year={2015},
  publisher={SIAM}
}

@article{radcliffe_vertex,
  title={Vertex isoperimetric inequalities for a family of graphs on $\mathbb{Z}^k$},
  author={Radcliffe, J and Veomett, E},
  journal={Electronic Journal of Combinatorics},
  volume={19},
  number={2},
  pages={P45},
  year={2012}
}

@article{ruzsa_asymptotics,
  title={Sets of sums and commutative graphs},
  author={Ruzsa, IZ},
  journal={Studia Scientiarum Mathematicarum Hungarica},
  volume={30},
  number={1},
  pages={127--148},
  year={1995},
  publisher={Budapest: Akademiai Kiado, 1966-}
}

@article{strachan2025edge,
  title={Edge isoperimetry of lattices},
  author={Strachan, Cameron and Swanepoel, Konrad},
  journal={arXiv preprint arXiv:2503.09591},
  year={2025}
}

@article{tsukerman_conjecture,
  title={A General Method to Determine Limiting Optimal Shapes for Edge-Isoperimetric Inequalities},
  author={Veomett, Ellen and Tsukerman, Emmanuel},
  journal={The Electronic Journal of Combinatorics},
  pages={P1--26},
  year={2017}
}

@article{wang_l1,
  title={Discrete isoperimetric problems},
  author={Wang, Da-Lun and Wang, Ping},
  journal={SIAM Journal on Applied Mathematics},
  volume={32},
  number={4},
  pages={860--870},
  year={1977},
  publisher={SIAM}
}

@book{zhao_book,
  title={Graph Theory and Additive Combinatorics: Exploring Structure and Randomness},
  author={Zhao, Yufei},
  year={2023},
  publisher={Cambridge University Press}
}

\end{document}